\documentclass[12pt,oneside,english]{amsart}
\usepackage[T1]{fontenc}
\usepackage[latin1]{inputenc}
\usepackage{amssymb}

\makeatletter

\providecommand{\LyX}{L\kern-.1667em\lower.25em\hbox{Y}\kern-.125emX\@}

 \theoremstyle{plain}    
 \newtheorem{thm}{Theorem} 
 \newcommand{\lyxaddress}[1]{
   \par {\raggedright #1 
   \vspace{1.4em}
   \noindent\par}
 }

\newcommand{\Tr}{\operatorname{Tr}}
\newcommand{\id}{\operatorname{id}}
\makeatother
\usepackage{babel}
\begin{document}

\title{On multiplicity and free absorption for free Araki-Woods factors.}

\date{\today}

\lyxaddress{\address{Department of Mathematics, UCLA, Los Angeles, CA 90095}}

\email{shlyakht@math.ucla.edu}

\thanks{This material is based upon work performed for the Clay Mathematics
Institute.}

\author{Dimitri Shlyakhtenko}

\begin{abstract}
We show that Ozawa's recent results on solid von Neumann algebras
imply that there are free Araki-Woods factors, which fail to have
free absorption. We also show that a free Araki-Woods factors $\Gamma (\mu ,n)$
associated to a measure and a multiplicity function $n$ may non-trivially
depend on the multiplicity function. 
\end{abstract}
\maketitle
In \cite{shlyakht:quasifree:big} we have associated to the absolute
continuity class of a symmetric measure $\mu $ on $\mathbb{R}$,
$\mu (X)=\mu (-X)$, and a multiplicity function $n:\mathbb{R}\to \mathbb{N}\cup \{+\infty \}$,
$n(x)=n(-x)$, a von Neumann algebra denoted by $\Gamma (\mu ,n)$%
\footnote{Our original notation for this algebra is $\Gamma (\mathcal{H}_{\mathbb{R}},U_{t})''$,
where $t\mapsto U_{t}$ is an orthogonal representation of $\mathbb{R}$
on a real Hilbert space. This data is of course completely equivalent
to specifying the spectral measure $\mu $ and multiplicity function
$n$ of such a representation. It is more convenient to use the notation
$\Gamma (\mu ,n)$ in the present paper.%
}. This construction is the free probability analog of the Araki-Woods
construction of hyperfinite type III factors using quasi-free states
\cite{araki,araki-woods}. The algebras $\Gamma (\mu ,n)$ are called
free Araki-Woods factors; they are factors of type III if $\mu $
is not concentrated on the set $\{0\}$. The complete classification
of free Araki-Woods factors remains an interesting outstanding problem,
and was achieved in \cite{shlyakht:quasifree:big} only for atomic
measures $\mu $. In our recent work \cite{shlyakht:fullfactor} we
have introduced new invariants for type III factors, which show that
the algebras $\Gamma (\mu ,n)$ depend to a large extent on the absolute
continuity class of the measure $\mu $. On the other hand, it remained
completely unknown in general whether the multiplicity function $n$
matters at all. In the case that $\mu $ is atomic and not concentrated
on $\{0\}$, we have shown in \cite{shlyakht:quasifree:big} that
$\Gamma (\mu ,n)$ is independent of the multiplicity function, leading
one to wonder if the same holds in general. 

Our construction of a free Araki-Woods factor $\Gamma (\mu ,n)$ produces
also a state $\phi _{\mu ,n}:\Gamma (\mu ,n)\to \mathbb{C}$, which
is the analog of the quasi-free state in the classical Araki-Woods
construction. We have shown in \cite{shlyakht:quasifree:big} that
the reduced free product of free Araki-Woods factors is given by:\[
(\Gamma (\mu ,n),\phi _{\mu ,n})*(\Gamma (\mu ',n'),\phi _{\mu ',n'})\cong (\Gamma (\mu +\mu ',n+n'),\phi _{\mu +\mu ',n+n'}),\]
with the isomorphism preserving the indicated states. 

Independence of $n$ then says in particular that\[
(\Gamma (\mu ,n),\phi _{\mu ,n})*(\Gamma (\mu ,n),\phi _{\mu ,n})\cong \Gamma (\mu ,n)\]
(via an isomorphism that does not necessarily take the free product
state to $\phi _{\mu ,n}$). 

Another property of free Araki-Woods factors which holds for $\mu $
atomic and not concentrated on $\{0\}$ is the {}``free absorption''
property:\[
\Gamma (\mu ,n)\cong (\Gamma (\mu ,n),\phi _{\mu ,n})*(L(\mathbb{F}_{\infty }),\tau ),\]
where $L(\mathbb{F}_{\infty })$ denotes the von Neumann algebra of
a free group on an infinite number of generators.

In this note we rely on the amazing recent result of Ozawa \cite{ozawa:solid}
showing that von Neumann algebras of hyperbolic groups (in particular,
of free groups) are {}``solid'': the relative commutant of any diffuse
unital subalgebra is hyperfinite. From this result, we derive the
following theorem, which shows that in the case that $\mu $ is non-atomic,
free absorption and independence of the multiplicity function need
not hold.

\begin{thm}
Let $\lambda $ denote the Lebesgue measure on $\mathbb{R}$, and
let $\delta _{0}$ denote the delta measure at $0$. Then:\\
(a) $\Gamma (\lambda +\delta _{0},1)\not \cong \Gamma (\lambda +\delta _{0},2)$.
Thus $\Gamma (\lambda +\delta _{0},n)$ depends on $n$.\\
(b) $\Gamma (\lambda ,1)\not \cong (\Gamma (\lambda ,1),\phi _{\lambda ,1})*(L(\mathbb{F}_{\infty }),\tau )$.
Thus free absorption fails.
\end{thm}
\begin{proof}
It is known that the core of $\Gamma (\lambda ,1)$ is isomorphic
to $L(\mathbb{F}_{\infty })\otimes B(\ell ^{2})$ (see \cite{shlyakht:amalg}).
We'll first show that the core of $\Gamma (\lambda +\delta _{0},1)$
is isomorphic to $L(\mathbb{F}_{\infty })\otimes B(\ell ^{2})$. 

We have that $(\Gamma (\lambda +\delta _{0},1),\phi _{\lambda +\delta _{0},1})\cong (\Gamma (\lambda ,1),\phi _{\lambda ,1})*(\Gamma (\delta _{0},1),\phi _{\delta _{0},1})$.
Thus by \cite{shlyakht:amalg}, the core of $\Gamma (\lambda +\delta _{0},1)$
is isomorphic to the amalgamated free product\[
(\Gamma (\lambda ,1)\rtimes _{\sigma ^{\phi _{\lambda ,1}}}\mathbb{R},E_{1})*_{L(\mathbb{R})}(\Gamma (\delta _{0},1)\rtimes _{\sigma ^{\phi _{\delta _{0},1}}}\mathbb{R},E_{2}).\]
Here $E_{j}$ are the conditional expectations onto the group von
Neumann algebra of $\mathbb{R}$ (viewed as a subalgebra of the crossed
product), associated to the states $\phi _{\lambda ,1}$ and $\phi _{\delta _{0},1}$,
respectively.

Since $\Gamma (\delta _{0},1)\cong L^{\infty }[0,1]$ by \cite{shlyakht:quasifree:big},
we have that the core of $\Gamma (\lambda +\delta _{0},1)$ is isomorphic
to\[
(\Gamma (\lambda ,1)\rtimes _{\sigma ^{\phi _{\lambda ,1}}}\mathbb{R},E_{1})*_{L(\mathbb{R})}(L(\mathbb{R})\otimes A,\id \otimes \tau ),\]
where $A\cong L^{\infty }[0,1]$, and $\tau $ is some normal faithful
trace on $A$.

Let $\Tr :L(\mathbb{R})\to \mathbb{C}$ be the semifinite trace associated
to integration against Lebesgue measure in the identification $L(\mathbb{R})\cong L^{\infty }(\mathbb{R})$.
Then by \cite{shlyakht:semicirc,shlyakht-popa:universal}, we have
that $(\Gamma (\lambda ,1)\rtimes \mathbb{R},E_{1})$ is isomorphic
to $(\Phi (L(\mathbb{R}),\Tr ),E_{L(\mathbb{R})})$ (see \cite{shlyakht-popa:universal}
for this notation). In other words, the core of $\Gamma (\lambda ,1)$,
when viewed as a non-commutative probability space over $L(\mathbb{R})$
is generated by $L(\mathbb{R})$ and an operator-valued semicircular
variable $X$ of covariance $\eta :L^{\infty }(\mathbb{R})\to L^{\infty }(\mathbb{R})$,
where $\eta $ is a completely-positive map given by\[
\eta (f)(x)=\int f(y)\omega (x,y)dy,\]
where $\omega >0$ is a certain function on $\mathbb{R}^{2}$ (i.e.,
$\omega (x,y)dxdy$ and $dxdy$ are mutually absolutely continuous).

Since the core of $\Gamma (\lambda +\delta _{0},1)$ is given as an
amalgamated free product, it follows that $X$ is free with amalgamation
over $L(\mathbb{R})$ from $A\otimes L(\mathbb{R})$. Thus by \cite{shlyakht:amalg},
if we view $X$ as an operator-valued random variable with values
in $A\otimes L(\mathbb{R})$, it is still semicircular and has covariance
$\eta \circ (\id \otimes \tau )\sim \Tr \circ (\id \otimes \tau ):L(\mathbb{R})\otimes A\to \mathbb{C}$.
But $(L(\mathbb{R})\otimes A,\Tr \otimes (\id \otimes \tau ))\cong (L(\mathbb{R}),\Tr )$.
Thus the core of $\Gamma (\lambda +\delta _{0},1)$ is isomorphic
to $\Phi (L(\mathbb{R}),\Tr )\cong L(\mathbb{F}_{\infty })\otimes B(\ell ^{2})$
by \cite{shlyakht:amalg}.

We now prove the non-isomorphisms claimed in (a) and (b). It is sufficient
to show that the cores of $\Gamma (\lambda +\delta _{0},2)$ and $(\Gamma (\lambda ,1),\phi _{\lambda ,1})*(L(\mathbb{F}_{\infty }),\tau )$
are not isomorphic to $L(\mathbb{F}_{\infty })\otimes B(\ell ^{2})$.

The core of $\Gamma (\lambda +\delta _{0},2)$ is isomorphic, as before,
to\[
(\Gamma (\lambda ,2)\rtimes _{\sigma ^{\phi _{\lambda ,2}}}\mathbb{R},E_{3})*_{L(\mathbb{R})}(\Gamma (\delta _{0},2)\rtimes _{\sigma ^{\phi _{\delta _{0},2}}}\mathbb{R},E_{4}).\]
But by \cite{shlyakht:quasifree:big}, $(\Gamma (\delta _{0},2),\phi _{\delta _{0},2})\cong (L(\mathbb{F}_{2}),\tau )$,
and the core of $\Gamma (\lambda +\delta _{0},2)$ is given by\[
(\Gamma (\lambda ,2)\rtimes _{\sigma ^{\phi _{\lambda ,2}}}\mathbb{R},E_{3})*_{L(\mathbb{R})}(L(\mathbb{F}_{2})\otimes L(\mathbb{R}),\tau \otimes \id ).\]
Now choose an abelian diffuse subalgebra $B\subset L(\mathbb{R})$,
so that the restriction of the trace on the core of $\Gamma (\lambda +\delta _{0},2)$
to $B$ is finite. Denote by $p$ the unit of $B$. Then the relative
commutant of $B=pBp$ inside of $p(\Gamma (\lambda ,2)\rtimes _{\sigma ^{\phi _{\lambda ,2}}}\mathbb{R},E_{3})*_{L(\mathbb{R})}(L(\mathbb{F}_{2})\otimes L(\mathbb{R}),\tau \otimes \id )p$
contains $L(\mathbb{F}_{2})$ and so is not hyperfinite. Thus $p(\Gamma (\lambda ,2)\rtimes _{\sigma ^{\phi _{\lambda ,2}}}\mathbb{R},E_{3})*_{L(\mathbb{R})}(L(\mathbb{F}_{2})\otimes L(\mathbb{R}),\tau \otimes \id )p$
is not solid. Thus the core of $\Gamma (\lambda +\delta _{0},2)$
is not a free group factor, by Ozawa's result \cite{ozawa:solid}.
This proves (a).

The core of $(\Gamma (\lambda ,1),\phi _{\lambda ,1})*(L(\mathbb{F}_{\infty },\tau ))$
is, as before, given by\[
(\Gamma (\lambda ,1)\rtimes _{\sigma ^{\phi _{\lambda ,1}}}\mathbb{R},E_{5})*_{L(\mathbb{R})}(L(\mathbb{F}_{\infty })\otimes L(\mathbb{R}),\tau \otimes \id ).\]
Using exactly the same argument as before, we find that the core is
not $L(\mathbb{F}_{\infty })\otimes B(\ell ^{2})$.
\end{proof}
Note that the core of $\Gamma (\lambda +\delta _{0},2)$ is related
to the factor $L(\mathbb{F}_{\infty })*(L(\mathbb{F}_{2})\otimes L^{\infty }[0,1])$
(the core is, in some loose sense, {}``$(L(\mathbb{F}_{\infty })\otimes B(\ell ^{2}))*(L(\mathbb{F}_{2})\otimes L^{\infty }(\mathbb{R}))$'').
This explains our interest in this factor (see the remark in \cite{ozawa:solid}). 

There is, however, a situation where free absorbtion holds, and the
mutiplicity function does not matter. For a measure $\nu $ and a
number $\lambda \in \mathbb{R}$, write $\nu _{\lambda }$ for the
translate of $\nu $ by $\lambda $, and $\nu ^{\lambda }$ for $\nu _{\lambda }+\nu _{-\lambda }$. 

\begin{thm}
Assume that $\mu =\mu _{c}+\mu _{a}$, where $\mu _{a}$ is atomic
and $\mu _{c}$ has no atoms. Assume that $\mu _{a}\neq 0$, and let
$\Gamma =\{\gamma _{1},\gamma _{2},\ldots \}$ be the subgroup of
$\mathbb{R}$ generated by its support. Assume that $\Gamma \neq \{0\}$.
Let $\nu =\sum _{j}2^{-j}\mu ^{\gamma _{j}}$. Then there is a state-preserving
isomorphism \[
(\Gamma (\mu ,n),\phi _{\mu ,n})\cong (\Gamma (\nu ,\infty ),\phi _{\nu ,\infty }).\]
In particular, $\Gamma (\mu ,n)$ does not depend on $n$, and free
absorbtion holds for $(\Gamma (\mu ,n),\phi _{\mu ,n})$. 
\end{thm}
\begin{proof}
We first note that the matricial model argument in \cite{shlyakht:quasifree:big}
implies that if $\phi $ is a normal state on $B(\ell ^{2})$ given
by\[
\theta (T)=\Tr (DT),\]
where $D$ is a trace-class diagonal matrix with eigenvalues $\exp (-\lambda _{1}),\exp (-\lambda _{2}),\ldots $,
with $\lambda _{1},\lambda _{2},\ldots $ non-negative, then for any
$\nu $ and $m$,\[
(B(\ell ^{2}),\theta )*(\Gamma (\nu ,m),\phi _{\nu ,m})\cong (\Gamma (\sum 2^{-j}\nu ^{\lambda _{j}},\infty )\otimes B(\ell ^{2}),\phi _{\sum 2^{-j}\nu ^{\lambda _{j}},\infty }\otimes \theta ).\]
Moreover, the isomorphism is state-preserving.

It follows that if the set $\{0\}\cup \{\pm \lambda _{j}:j=1,2,\ldots \}$
is an additive group, then\begin{eqnarray*}
(B(\ell ^{2}),\theta )*(\Gamma (\sum 2^{-j}\nu ^{\lambda _{j}},\infty ) & \cong  & (\Gamma (\sum 2^{-j}\nu ^{\lambda _{j}},\infty )\otimes B(\ell ^{2}),\phi _{\sum 2^{-j}\nu ^{\lambda _{j}},\infty }\otimes \theta )\\
 & \cong  & (B(\ell ^{2}),\theta )*(\Gamma (\nu ,m),\phi _{\nu ,m}),
\end{eqnarray*}
via state-preserving isomorphisms.

Returning now to our situation, note that $(\Gamma (\mu ,n),\phi _{\mu ,n})\cong (\Gamma (\mu _{c},n),\phi _{\mu _{c},n})*(\Gamma (\mu _{a},n),\phi _{\mu _{a},n})$.
By \cite{shlyakht:quasifree:big}, $(\Gamma (\mu _{a},n),\phi (\mu _{a},n))\cong (\Gamma (\mu _{a},n),\phi (\mu _{a},n))*(B(\ell ^{2}),\theta )$,
where $\theta (T)=\Tr (DT)$, and $D$ is a diagonal operator with
eigenvalues $\{e^{\gamma _{j}}:\gamma _{j}<0\}$. Thus $(\Gamma (\mu ,n),\phi _{\mu ,n})\cong (\Gamma (\mu ,n),\phi _{\mu ,n})*(B(\ell ^{2}),\theta )$.
By the argument above, it follows that\begin{eqnarray*}
(\Gamma (\mu ,n),\phi _{\mu ,n}) & \cong  & (\Gamma (\mu ,n),\phi _{\mu ,n})*(B(\ell ^{2}),\theta )\\
 & \cong  & (\Gamma (\sum 2^{-j}\mu ^{\gamma _{j}},\infty ),\phi _{\sum 2^{-j}\mu ^{\gamma _{j}},\infty })*(B(\ell ^{2}),\theta )\\
 & \cong  & (\Gamma (\sum 2^{-j}\mu ^{\gamma _{j}},\infty ),\phi _{\sum 2^{-j}\mu ^{\gamma _{j}},\infty }).
\end{eqnarray*}
Since $\sum 2^{-j}\mu ^{\gamma _{j}}$ has an atom at zero, it follows
that\begin{eqnarray*}
(\Gamma (\mu ,n),\phi _{\mu ,n}) & \cong  & (\Gamma (\sum 2^{-j}\mu ^{\gamma _{j}},\infty ),\phi _{\sum 2^{-j}\mu ^{\gamma _{j}},\infty })\\
 & \cong  & (\Gamma (\sum 2^{-j}\mu ^{\gamma _{j}},\infty ),\phi _{\sum 2^{-j}\mu ^{\gamma _{j}},\infty })*(L(\mathbb{F}_{\infty }),\tau )\\
 & \cong  & (\Gamma (\mu ,n),\phi _{\mu ,n})*(L(\mathbb{F}_{\infty }),\tau ),
\end{eqnarray*}
so that free absorption holds.
\end{proof}
\bibliographystyle{amsplain}

\providecommand{\bysame}{\leavevmode\hbox to3em{\hrulefill}\thinspace}
\providecommand{\MR}{\relax\ifhmode\unskip\space\fi MR }
\providecommand{\MRhref}[2]{%
  \href{http://www.ams.org/mathscinet-getitem?mr=#1}{#2}
}
\providecommand{\href}[2]{#2}

\end{document}